# General class of mixture of some densities

**Salah Hamza Abid**

Department of Mathematics, College of Education, University of Mustansiriyah, Baghdad, Iraq
Corresponding author; e-mail: abidsalah@gmail.com

**Abstract**

In this paper, a general class of mixture of some densities is proposed. The proposed class contains some of classical and weighted distributions as special cases. Formulas for each of cumulative distribution function, reliability function, hazard rate function, rth raw moments function, characteristic function, stress-strength reliability and Tsallis entropy of order $\alpha$ are derived.

______________________

**1. Introduction**

Recently, many researchers have investigated for new probability distributions. Various methods have been proposed for deriving these distributions. All of these efforts are to find flexible distributions that are more representative than others for set(s) of data. Generalized distributions are often the ones that do so.

Mixture models are very important family of probability distributions. Cluster analysis and classification analysis and modeling failure data are most important uses of mixture models. In addition to their use in other statistical aspects due to its high flexibility. The important step for active use of mixture model is the correct determining of the number of subpopulations that the data supports and their distributions. In many practical applications the existing data can be viewed as arising from various populations. It is important to formulate a model that represents all subpopulations. This model will admits us to treat with the samples from different subpopulations.

In other words, to represent a population that is have subpopulations by a statistical model, we used usually the mixture densities, where the mixture parts are the distributions on the subpopulations and the weights are the ratios of each subpopulation in the overall population.

The probability density function (pdf) of the proposed general class of mixture of some densities is,

$$f(x) = c \left(\sum_{i=0}^{p} \theta_i x^i\right) e^{-\beta x^d}, \quad x > 0 \tag{1}$$

Where, $c$ is a constant which is determined as follows,

Since, $\int_0^\infty f(x)dx = c \sum_{i=0}^{p} \theta_i \int_0^\infty x^i e^{-\beta x^d} dx = 1$, by using the transformation $y = x^d$, then,

$$\int_0^\infty x^i e^{-\beta x^d} dx = (1/d) \int_0^\infty y^{(i+1)/d-1} e^{-\beta y} dy = \frac{\Gamma((i+1)/d)}{d\,\beta^{(i+1)/d}}\ , \text{And then,}$$

$$c = \frac{d}{\sum_{i=0}^p \theta_i \frac{\Gamma((i+1)/d)}{\beta^{(i+1)/d}}}\ , \tag{2}$$

The mixing proportions can be determined as follows,

$$MP_s = \frac{\theta_s \frac{\Gamma((s+1)/d)}{\beta^{(s+1)/d}}}{\sum_{i=0}^p \theta_i \frac{\Gamma((i+1)/d)}{\beta^{(i+1)/d}}} \tag{3}$$

So, the cumulative distribution function (cdf) will be,

$$F(x) = c \sum_{i=0}^p \theta_i \int_0^x x^i e^{-\beta x^d} dx\ , \text{by using the transformation } y = \beta x^d, \text{then,}$$

$$\int_0^x x^i e^{-\beta x^d} dx = \frac{1}{d\,\beta^{(i+1)/d}} \int_0^{\beta x^d} y^{(i+1)/d-1} e^{-y} dy$$

$$= \frac{1}{d\,\beta^{(i+1)/d}} \gamma\big((i+1)/d,\ \beta x^d\big)\ , \text{and then,}$$

$$F(x) = \frac{\sum_{i=0}^p \theta_i \beta^{-(i+1)/d} \gamma\big((i+1)/d, \beta x^d\big)}{\sum_{j=0}^p \theta_j \beta^{-(j+1)/d} \Gamma((j+1)/d)} \tag{4}$$

There are a lot sub-models of the proposed class. Most of them are contained in table (1) in details.

Also the rth raw moment for the class is as in the following theorem,

**Theorem 1** the rth raw moment of class in (1) is,

$$E(X^r) = \frac{\sum_{i=0}^p \theta_i \beta^{-(i+r+1)/d} \Gamma((i+r+1)/d)}{\sum_{j=0}^p \theta_j \beta^{-(j+1)/d} \Gamma((j+1)/d)} \tag{5}$$

**Proof:**

Since, $E(X^r) = \frac{d}{\sum_{j=0}^p \theta_j \frac{\Gamma((j+1)/d)}{\beta^{(j+1)/d}}} \sum_{i=0}^p \theta_i \int_0^\infty x^{r+i} e^{-\beta x^d} dx$, then, by using the transformation $y = x^d$, we get $E(X^r) = \frac{1}{\sum_{j=0}^p \theta_j \frac{\Gamma((j+1)/d)}{\beta^{(j+1)/d}}} \sum_{i=0}^p \theta_i \int_0^\infty y^{(r+i+1)/d-1} e^{-\beta y} dy$ , and then, the proof is complete.

**Corollary 1** the characteristic function of class in (1) is,

$$\varphi_X(t) = \sum_{r=0}^\infty \frac{it^r}{r!} \frac{\sum_{i=0}^p \theta_i \beta^{-(i+r+1)/d} \Gamma\big((i+r+1)/d, \beta X^d\big)}{\sum_{j=0}^p \theta_j \beta^{-(j+1)/d} \Gamma((j+1)/d)} \tag{6}$$

**Proof:** It is clear by the fact that, $\varphi_X(t) = E(e^{itX}) = \sum_{r=0}^\infty \frac{it^r}{r!} E(X^r)$.

**Table 1** Some sub-models of proposed class

| Distribution | Parameters | Author(s) |
|---|---|---|
| $Exp(\beta)$ | $\theta_0 = 1, \theta_i = 0$ for $(i = 1,2,...,p)$ and $d = 1$. | common |
| $Ga(\alpha, \beta)$, $\alpha$ is an integer | $p = \alpha - 1, \theta_i = 0$ for $(i = 0,1,...,\alpha - 2)$, $\theta_{\alpha-1} = 1$ and $d = 1$. | common |
| $We(d, \beta)$, $d$ is an integer | $p = d - 1, \theta_i = 0$ for $(i = 0,1,...,d - 2), \theta_{d-1} = 1$. | Weibull 1951 [23] Frechet 1927 [1] |
| $GGa(d, \beta, b)$, $b$ is an integer | $p = b - 1, \theta_i = 0$ for $(i = 0,1,...,b - 2), \theta_{b-1} = 1$. | Stacy 1962 [22] |
| $Lnd(\beta)$ {Mixture of $Exp(\beta)$ and $Ga(2,\beta)$} with $MP = \beta/(\beta + 1)$ | $d = 1, p = 1, \theta_0 = 1, \theta_1 = 1$. | Lindley 1958 [2] |
| $Aradhana(\beta)$ {Mixture of $Exp(\beta), Ga(2,\beta)$ and $Ga(3,\beta)$} with $MP_1 = \beta^2/(\beta^2 + 2\beta + 2)$ and $MP_2 = 2\beta/(\beta^2 + 2\beta + 2)$ | $d = 1, p = 2, \theta_0 = 1, \theta_1 = 2, \theta_2 = 1$. | Shanker (2016) [9] |
| $Ishita(\beta)$ {Mixture of $Exp(\beta)$ and $Ga(3,\beta)$} with $MP = \beta^3/(\beta^3 + 2)$ | $d = 1, p = 2, \theta_0 = \beta, \theta_1 = 0, \theta_2 = 1$. | Shanker and Shukla, (2017) [17] |
| $Akash(\beta)$ {Mixture of $Exp(\beta)$ and $Ga(3,\beta)$} with $MP = \beta^2/(\beta^2 + 2)$ | $d = 1, p = 2, \theta_0 = 1, \theta_1 = 0, \theta_2 = 1$. | Shanker (2015) [7] |
| $Amarendra(\beta)$ {Mixture of $Exp(\beta), Ga(2,\beta), Ga(3,\beta)$ and $Ga(4,\beta)$} with $MP_1 = \beta^3/(\beta^3 + \beta^2 + 2\beta + 6)$ and $MP_2 = \beta^2/(\beta^3 + \beta^2 + 2\beta + 6)$ and $MP_3 = 2\beta/(\beta^3 + \beta^2 + 2\beta + 6)$. | $d = 1, p = 3, \theta_i = 1$ for $(i = 0,1,2,3)$ | Shanker (2016) [11] |
| $Sujatha(\beta)$ {Mixture of $Exp(\beta), Ga(2,\beta)$ and $Ga(3,\beta)$} with $MP_1 = \beta^2/(\beta^2 + \beta + 2)$ and $MP_2 = \beta/(\beta^2 + \beta + 2)$. | $d = 1, p = 2, \theta_i = 1$ for $(i = 0,1,2)$ | Shanker (2016) [10] |
| $Shanker(\beta)$ {Mixture of $Exp(\beta)$ and $Ga(2,\beta)$} with $MP = \beta^2/(\beta^2 + 1)$ | $d = 1, p = 1, \theta_0 = \beta, \theta_1 = 1$. | Shanker (2015) [8] |
| $Akshaya(\beta)$ {Mixture of $Exp(\beta), Ga(2,\beta), Ga(3,\beta)$ and $Ga(4,\beta)$} with $MP_1 = \beta^3/(\beta^3 + 3\beta^2 + 6\beta + 6)$ and $MP_2 = 3\beta^2/(\beta^3 + 3\beta^2 + 6\beta + 6)$ and $MP_3 = 6\beta/(\beta^3 + 3\beta^2 + 6\beta + 6)$. | $d = 1, p = 3, \theta_0 = \theta_3 = 1, \theta_1 = \theta_2 = 3$. | Shanker (2017) [14] |
| $Suja(\beta)$ {Mixture of $Exp(\beta)$ and $Ga(5,\beta)$} with $MP = \beta^4/(\beta^4 + 24)$ | $d = 1, p = 4, \theta_i = 1$ for $(i = 0,4)$ and $\theta_j = 0$ for $(j = 1,2,3)$ | Shanker (2017) [19] |
| $Devyaa(\beta)$ {Mixture of $Exp(\beta), Ga(2,\beta), Ga(3,\beta), Ga(4,\beta)$ and $Ga(5,\beta)$} with $MP_1 = \beta^4/(\beta^4 + \beta^3 + 2\beta^2 + 6\beta + 24)$ and $MP_2 = \beta^3/(\beta^4 + \beta^3 + 2\beta^2 + 6\beta + 24)$ and | $d = 1, p = 4, \theta_i = 1$ for $(i = 0,1,2,3,4)$ | Shanker (2016) [13] |
| $MP_3 = 2\beta^2/(\beta^4 + \beta^3 + 2\beta^2 + 6\beta + 24)$ and $MP_4 = 6\beta/(\beta^4 + \beta^3 + 2\beta^2 + 6\beta + 24)$. | | |
| $QLnd(\alpha,\beta)$ {Mixture of $Exp(\beta)$ and $Ga(2,\beta)$} with $MP = \alpha/(\alpha + 1)$ | $d = 1, p = 1, \theta_0 = \alpha, \theta_1 = \beta$. | Shanker, and Mishra, (2013) [4] |
| $GSujatha(\alpha,\beta)$ {Mixture of $Exp(\beta)$, $Ga(2,\beta)$ and $Ga(3,\beta)$} with $MP_1 = \beta^2/(\beta^2 + \beta + 2\alpha)$ and $MP_2 = \beta/(\beta^2 + \beta + 2\alpha)$. | $d = 1, p = 2, \theta_i = 1$ for $(i = 0,1), \theta_2 = \alpha$. | Shanker (2017) [15] |
| $Rani(\beta)$ {Mixture of $Exp(\beta)$ and $Ga(5,\beta)$} with $MP = \beta^5/(\beta^5 + 24)$ | $d = 1, p = 5, \theta_0 = \beta, \theta_4 = 1, \theta_i = 0$ for $(i = 1,2,3)$. | Shanker (2017) [16] |
| $Garima(\beta)$ {Mixture of $Exp(\beta)$ and $Ga(2,\beta)$} with $MP = (\beta + 1)/(\beta + 2)$ | $d = 1, p = 1, \theta_0 = \beta + 1, \theta_1 = \beta$. | Shanker (2016) [12] |
| $Janardan(\alpha,\mu)$ {Mixture of $Exp(\mu/\alpha)$ and $Ga(2,\mu/\alpha)$} with $MP = \mu/(\mu + \alpha^2)$ | $d = 1, \beta = \mu/\alpha, p = 1, \theta_0 = 1, \theta_1 = \alpha$. | Shanker (2013) [5] |
| $Om(\beta)$ {Mixture of $Exp(\beta), Ga(2,\beta), Ga(3,\beta), Ga(4,\beta)$ and $Ga(5,\beta)$} with $MP_1 = \beta^4/(\beta^4 + 4\beta^3 + 12\beta^2 + 24\beta + 24)$ and $MP_2 = 4\beta^3/(\beta^4 + 4\beta^3 + 12\beta^2 + 24\beta + 24)$ and $MP_3 = 12\beta^2(\beta^4 + 4\beta^3 + 12\beta^2 + 24\beta + 24)$ and $MP_4 = 24\beta/(\beta^4 + 4\beta^3 + 12\beta^2 + 24\beta + 24)$. | $d = 1, p = 4, \theta_i = 1$ for $(i = 0,4), \theta_j = 4$ for $(i = 1,3)$ and $\theta_2 = 6$ | Shanker and Shukla, (2018) [20] |
| $Sushila(\alpha,\mu)$ {Mixture of $Exp(\mu/\alpha)$ and $Ga(2,\mu/\alpha)$} with $MP = \mu/(\mu + 1)$ | $d = 1, \beta = \mu/\alpha, p = 1, \theta_0 = 1, \theta_1 = 1/\alpha$. | Shanker et al. (2013) [6] |
| $GAradhana(\alpha,\beta)$ {Mixture of $Exp(\beta), Ga(2,\beta)$ and $Ga(3,\beta)$} with $MP_1 = \beta^2/(\beta^2 + 2\alpha\beta + 2\alpha^2)$ and $MP_2 = 2\alpha\beta/(\beta^2 + 2\alpha\beta + 2\alpha^2)$ | $d = 1, p = 2, \theta_0 = 1, \theta_1 = 2\alpha, \theta_2 = \alpha^2$. | Welday and Shanker, (2018) [21] |
| $GLnd(\alpha,\lambda,\beta)$ {Mixture of $Exp(\beta)$ and $Ga(2,\beta)$} with $MP = \alpha\beta/(\alpha\beta + \lambda)$ | $d = 1, p = 1, \theta_0 = \alpha, \theta_1 = \lambda$. | Shanker et al. (2017) [18] |
| $XGa(\beta)$ {Mixture of $Exp(\beta)$ and $Ga(3,\beta)$} with $MP = \beta/(\beta + 1)$ | $d = 1, p = 2, \theta_0 = 1, \theta_1 = 0, \theta_2 = \beta/2$. | Sen et al. (2016) [3] |

## 2. Stress-Strength model and Tsallis entropy

Let $X$ and $Y$ are the independent strength and stress random variables distributed as class in (1) with completely different parameters, then the Stress-Strength model for the class is as in the following theorem,

**Theorem 2** the Stress-Strength model of the class in (1) is,

$$R = \frac{c}{d\,\tau}\sum_{i=0}^{p}\sum_{j=0}^{p^*}\theta_i\theta_j^*\,\beta^{*-(j+1)/d^*}\sum_{k=0}^{\infty}\frac{(-1)^k \beta^{(j+1)/d^*+k}\Gamma((i+j+kd^*+2)/d)}{k!\,((j+1)/d^*+k)\,\beta^{(i+j+kd^*+2)/d}} \tag{7}$$

Where, $\tau = \left\{\sum_{j=0}^{p}\theta_j^*\,\beta^{*-(j+1)/d^*}\Gamma((j+1)/d^*)\right\}^{-1}$

**Proof:**

$R = P(Y<X) = \int_0^{\infty} f_X(x)\,F_Y(x)\,dx$

$= \int_0^{\infty}(c/\tau)\left(\sum_{i=0}^{p}\theta_i x^i\right)e^{-\beta x^d}\sum_{j=0}^{p^*}\theta_j^*\,\beta^{*-(j+1)/d^*}\gamma\big((j+1)/d^*,\,\beta^* x^{d^*}\big)\,dx$

$= \frac{c}{\tau}\sum_{i=0}^{p}\sum_{j=0}^{p^*}\theta_i\theta_j^*\,\beta^{*-(j+1)/d^*}\int_0^{\infty}x^i e^{-\beta x^d}\gamma\big((j+1)/d^*,\,\beta^* x^{d^*}\big)dx$

Since, $(u,z) = \sum_{k=0}^{\infty}\frac{(-1)^k z^{u+k}}{k!\,(u+k)}$ (8)

Then,

$\int_0^{\infty} x^i e^{-\beta x^d}\gamma\big((j+1)/d^*,\,\beta^* x^{d^*}\big)dx = \int_0^{\infty} x^i e^{-\beta x^d}\sum_{k=0}^{\infty}\frac{(-1)^k\left(\beta^* x^{d^*}\right)^{(j+1)/d^*+k}}{k!\,((j+1)/d^*+k)}dx$

$\qquad = \sum_{k=0}^{\infty}\frac{(-1)^k \beta^{*(j+1)/d^*+k}}{k!\,((j+1)/d^*+k)}\int_0^{\infty}x^{i+d^*k+j+1}e^{-\beta x^d}dx$

Now, since $\int_0^{\infty} x^{b-1}e^{-\beta x^d}dx = \frac{\Gamma(b/d)}{d\,\beta^{b/d}}$, then the proof is complete.

The Tsallis entropy of order $\alpha$ which is proposed by Tsallis (1988), $TE = \frac{1}{\alpha-1}\left(1 - \int_0^{\infty} f^{\alpha}(x)\,dx\right)$, where $\alpha$ is real number ($\alpha \geq 0,\,\alpha \neq 1$), is derived for the class is as in the following theorem,

**Theorem 3** if $\alpha$ is positive integer ($\alpha \geq 0,\,\alpha \neq 1$), the Tsallis entropy of the class in (1) is,

$$TE = \frac{1}{\alpha-1}\left(1 - \frac{c^{\alpha}}{d}\sum_{i_1=0}^{p}\cdots\sum_{i_\alpha=0}^{p}\theta_{i_1}\cdots\theta_{i_\alpha}\frac{\Gamma((i_1+\cdots+i_\alpha+1)/d)}{(\beta\alpha)^{(i_1+\cdots+i_\alpha+1)/d}}\right) \tag{8}$$

**Proof:**

Since,

$\int_0^{\infty} f^{\alpha}(x)\,dx = \int_0^{\infty} c^{\alpha}\left(\sum_{i=0}^{p}\theta_i x^i\right)^{\alpha} e^{-\alpha\beta x^d}dx$

$\qquad = \int_0^{\infty} c^{\alpha}\sum_{i_1}^{p}\cdots\sum_{i_\alpha}^{p}(\theta_{i_1}\cdots\theta_{i_\alpha})x^{i_1+\cdots+i_\alpha}e^{-\alpha\beta x^d}dx$

$\qquad = c^{\alpha}\sum_{i_1}^{p}\cdots\sum_{i_\alpha}^{p}(\theta_{i_1}\cdots\theta_{i_\alpha})\int_0^{\infty} x^{i_1+\cdots+i_\alpha}e^{-\alpha\beta x^d}dx$

$$= \frac{c^\alpha}{d} \sum_{i_1=0}^{p} \cdots \sum_{i_\alpha=0}^{p} \theta_{i_1} \cdots \theta_{i_\alpha} \frac{\Gamma((i_1+\cdots+i_\alpha+1)/d)}{(\beta\alpha)^{(i_1+\cdots+i_\alpha+1)/d}}$$

Then the proof is complete.

### 3. Examples

In this section we will discuss two examples about the proposed class,

Example (1):

Following discussion is for the mixture distribution of $GG(1,\beta,2)$ and $GG(3,\beta,2)$ random variables.

1. The pdf of the distribution is $f(x) = c\,(1+x^2)e^{-\beta x^2}, x > 0$, with $c = \frac{1}{\frac{\Gamma(1/2)}{2\beta^{1/2}}} + \frac{1}{\frac{\Gamma(3/2)}{2\beta^{3/2}}}$

   Then $f(x) = \frac{\beta^{3/2}}{(2\beta+1)\sqrt{\pi}} (1+x^2)e^{-\beta x^2}, x > 0$.

2. The Mixing proportions are, $MP_1 = 2\beta/(2\beta+1)$ and $MP_2 = 1/(2\beta+1)$.

3. The cdf of the distribution is $F(x) = \frac{\beta\,\gamma(1/2,\beta x^2) + \gamma(2/3,\beta x^2)}{2(2\beta+1)\sqrt{\pi}}$

4. The rth raw moment function is, $E(X^r) = \frac{\beta\Gamma((r+1)/2) + \Gamma((r+3)/2)}{2(2\beta+1)\sqrt{\pi\beta^r}}$

   So, the mean and the variance of the distribution are,
   $$\mu = E(X) = \frac{\beta+1}{2(2\beta+1)\sqrt{\pi\beta}} \quad \text{and} \quad \sigma^2 = \frac{\beta/2 + 3/4}{2\beta(2\beta+1)} - \frac{(\beta+1)^2}{4(2\beta+1)^2\,\beta\sqrt{\pi}}$$

   And one can easily specified the skewness and kurtosis by using their relations with the above mean and variance side by side with the following third and fourth raw moments,
   $$E(X^3) = \frac{\beta+2}{2(2\beta+1)\sqrt{\pi\beta^3}} \quad \text{and} \quad E(X^4) = \frac{3\beta + 15/2}{8(2\beta+1)\beta^2}$$

5. The stress-strength reliability model is,

$$R = \frac{\beta^{3/2}}{(2\beta+1)\sqrt{\pi}} (1+x^2)e^{-\beta x^2} \left(\frac{\rho\,\gamma(1/2,\rho x^2) + \gamma(2/3,\rho x^2)}{2(2\rho+1)\sqrt{\pi}}\right) dx$$

$$= \frac{\beta^{3/2}}{2\pi(2\rho+1)(2\beta+1)} \left(\rho \int_0^\infty (1+x^2)e^{-\beta x^2}\,\gamma(1/2,\rho x^2) + \int_0^\infty (1+x^2)e^{-\beta x^2}\,\gamma(2/3,\rho x^2)\right) dx$$

By using (8), we get,

$$R = \frac{\beta^{\frac{3}{2}}}{2\pi(2\rho+1)(2\beta+1)} \left(\rho \int_0^\infty (1+x^2)e^{-\beta x^2} \sum_{p=0}^{\infty} \frac{(-1)^p (\rho x^2)^{p+\frac{1}{2}}}{p!\left(p+\frac{1}{2}\right)} + \int_0^\infty (1+x^2)e^{-\beta x^2} \sum_{q=0}^{\infty} \frac{(-1)^q (\rho x^2)^{q+\frac{3}{2}}}{q!\left(q+\frac{3}{2}\right)}\right) dx$$

$$= \frac{\beta^{\frac{3}{2}}}{2\pi(2\rho+1)(2\beta+1)} \left(\rho \sum_{p=0}^{\infty} \frac{(-1)^p \rho^{p+\frac{1}{2}}}{p!\left(p+\frac{1}{2}\right)} \int_0^\infty (1+x^2)e^{-\beta x^2} x^{2p+1} + \sum_{q=0}^{\infty} \frac{(-1)^q \rho^{q+\frac{3}{2}}}{q!\left(q+\frac{3}{2}\right)} \int_0^\infty (1+x^2)e^{-\beta x^2} x^{2q+3}\right) dx$$

$$= \frac{1}{2\pi(2\rho+1)} \left(\rho \sum_{p=0}^{\infty} \frac{(-1)^p \rho^{p+\frac{1}{2}}}{p!\left(p+\frac{1}{2}\right)} \left(\frac{\beta\Gamma(p+1)+\Gamma(p+2)}{2(2\beta+1)\sqrt{\pi\beta^{2p+1}}}\right) + \sum_{q=0}^{\infty} \frac{(-1)^q \rho^{q+\frac{3}{2}}}{q!\left(q+\frac{3}{2}\right)} \left(\frac{\beta\Gamma(q+2)+\Gamma(q+3)}{2(2\beta+1)\sqrt{\pi\beta^{2q+3}}}\right)\right)$$

Example (2):

Following discussion is for the mixture distribution of $Exp(\beta), Ga(2,\beta), \ldots, Ga(p+1, \beta)$ random variables.

1. The pdf of the distribution is $f(x) = c\,(a + bx)^p e^{-\beta x}, x > 0$.
2. $d = 1, \theta_s = C_s^p a^{p-s} b^s, \quad s = 0,1,\ldots,p$.
3. $c = \dfrac{1}{\sum_{i=0}^{p} C_i^p a^{p-i} b^i \frac{\Gamma(i+1)}{\beta^{i+1}}} = \dfrac{\beta}{a^p}\left(\sum_{i=0}^{p} \mathcal{P}_i^p (b/a\beta)^i\right)^{-1}$, where $\mathcal{P}$ denote a permutation.
4. The mixing proportions are,
$$MP_m = \frac{\mathcal{P}_m^p (b/a\beta)^m}{\sum_{i=0}^{p} \mathcal{P}_i^p (b/a\beta)^i}, \quad m = 0,1,\ldots,p \ , \text{ where } \sum_{m=0}^{p} MP_m = 1.$$
5. The cdf of the distribution is,
$$F(x) = \frac{\sum_{i=0}^{p} C_i^p a^{p-i} b^i \beta^{-(i+1)} \gamma\big((i+1),\beta x\big)}{\sum_{j=0}^{p} C_j^p a^{p-j} b^j \beta^{-(j+1)} \Gamma\big((j+1)\big)}$$
6. The rth raw moment function is,
$$E(X^r) = \frac{\sum_{i=0}^{p} C_i^p a^{p-i} b^i \beta^{-(i+r+1)} \Gamma\big((i+r+1)\big)}{\sum_{j=0}^{p} C_j^p a^{p-j} b^j \beta^{-(j+1)} \Gamma\big((j+1)\big)}$$
7. The stress-strength reliability model is,
$$R = \frac{c}{\tau} \sum_{i=0}^{p} \sum_{j=0}^{p^*} C_i^p a^{p-i} b^i\, C_j^{p^*} (a^*)^{p^*-j} (b^*)^j\, \beta^{*-(j+1)/d^*} \sum_{k=0}^{\infty} \frac{(-1)^k \beta^{(j+1)/d^*+k} \Gamma\big((i+j+kd^*+2)\big)}{k!\,((j+1)/d^*+k)\,\beta^{(i+j+kd^*+2)}}$$
Where, $\tau = \left\{\sum_{j=0}^{p} C_j^{p^*} (a^*)^{p^*-j} (b^*)^j \beta^{*-(j+1)/d^*} \Gamma\big((j+1)/d^*\big)\right\}^{-1}$

## 4. Summary and conclusions

Due to the great role that data analysis plays in various scientific fields, and because of the value of probability distributions in that, we notice an unparalleled acceleration of various developments towards deriving new distributions to fit the actual data of phenomena more closely. In this endeavor comes our research in which we present a general class of mixture distributions, where the most important statistical properties of this class are covered in addition to showing that many of the probability distributions represent special cases of it. Furthermore, we supported our work by illustrative examples.